\documentclass[letterpaper, 10 pt, conference]{ieeeconf} 

\IEEEoverridecommandlockouts                              
\usepackage{cite}

\usepackage{amsmath,amssymb,amsfonts,amsthm}
\usepackage[utf8]{inputenc}
\usepackage[T1]{fontenc}
\usepackage{booktabs}
\usepackage{graphics} 
\usepackage{epsfig} 
\usepackage{mathptmx} 
\usepackage{times} 
\usepackage{amsmath} 
\usepackage{amssymb}  
\usepackage{float}
\usepackage{xcolor}
\usepackage[most]{tcolorbox}

\title{\LARGE \bf
Integral Formulation of QENDy for Robust Nonlinear System Identification
}

\author{
Nikhil Saran$^{1}$, 
Sushant Pokhriyal$^{2}$, 
Stefan Klus$^{3}$,
Rushikesh Kamalapurkar$^{4}$,
and Joel A. Rosenfeld$^{5}$%
\thanks{*This work was partially supported by the following grants: AFOSR Program Awards FA9550-21-1-0134 and FA9550-20-1-0127 and NSF awards ECCS-2027976 and 2027999. It was also supported by the DOE/NIH/NSF ``Collaborative Research in Computational Neuroscience (CRCNS): Decomposing Neural Dynamics." Any opinions, findings and conclusions or recommendations expressed in this material are those of the author(s) and do not necessarily reflect the views of the sponsoring agencies.}
\thanks{$^{1}$Nikhil Saran (aka Nikhil Nikhil) is with the Department of Mathematics and Statistics at the University of South Florida, Email: nikhil36@usf.edu}
\thanks{$^{2}$Sushant Pokhriyal is an Assistant Professor at the Institute of Engineering and Technology, JK Lakshmipat University, Jaipur, India, Email: sushant.pokhriyal@jklu.edu.in}%
\thanks{$^{3}$Stefan Klus is an Associate Professor in the School of Mathematical \& Computer Sciences at the Heriot--Watt University, Edinburgh, Email: S.Klus@hw.ac.uk}
\thanks{$^{4}$Rushikesh Kamalapurkar is an Associate Professor at the University of Florida, Email:rkamalapurkar@ufl.edu}
\thanks{$^{5}$Joel A. Rosenfeld is an Associate Professor in the Department of Mathematics and Statistics at the University of South Florida, Email: rosenfeldj@usf.edu}%
}

\begin{document}

\maketitle
\thispagestyle{empty}
\pagestyle{empty}

\begin{abstract}
This manuscript proposes an integral formulation of the newly defined quadratic embedding method for identifying nonlinear systems (QENDy). In the original algorithm, trajectory data points along with their time derivatives are used. Methods for calculating time derivatives make the algorithm sensitive to noise. Our integral formulation does not use the time derivatives. This results in a more robust method to learn the dynamics.
\end{abstract}

\section{INTRODUCTION}

Many machine learning and data science techniques seek to uncover the underlying structure of a complex nonlinear dynamical system. Using only trajectory data, they identify the patterns of evolution of states over time. One of the prominent data-driven algorithms is SINDy \cite{brunton2016sindy}.
It is a sparse regression method that utilizes trajectory data, along with its time derivatives, and a set of basis functions. It states that only a few of the basis functions are actually needed for approximating the governing equations of the system. The algorithm starts with a larger set of basis functions and uses $l1$~$regularization$ to promote sparsity. The accuracy of this algorithm depends upon the choice of the library of functions, and its reliance on numerical derivatives makes it highly sensitive to measurement noise. Together, these factors can limit its robustness when applied to real-world data.
\par
Recently, Stefan Klus and Joel-Pascal N'konzi introduced system identification using quadratic embeddings of nonlinear dynamics, also known as QENDy \cite{klus2025quadratic}. The core idea is to obtain a finite-dimensional quadratic model of the dynamics in a higher-dimensional space. This $quadratization$ or $ polynomialization$ of the dynamics is always possible, provided the system comprises elementary functions such as polynomials, trigonometric, exponential, and rational functions \cite{kerner1981universal}, see also \cite{gu2011qlmor}. This representation is achieved through successively defining new variables and using the chain rule. For a simple example, take the differential equation 
\[ {\dot{x}} = x^3 \]
Introduce observables
\[
z_1 = x, \quad z_2 = x^2.
\]
Then, by the chain rule,
\[
\dot{z}_1 = \dot{x} = x^3 = z_1 z_2, \qquad
\dot{z}_2 = 2x \dot{x} = 2x^4 = 2 z_2^2.
\]
Thus, the system is embedded in the quadratic system
\[
\begin{split}
\dot{z}_1 &= z_1 z_2, \\
\dot{z}_2 &= 2z_{2}^{2}. \\
\end{split}
\]

Quadratic embedding is crucial because it transforms complex nonlinear dynamics into a structured quadratic form that is tractable for learning. It acts as a compromise between the linear but infinite-dimensional Koopman representation and the original finite-dimensional but fully nonlinear system \cite{colbrook2023multiverse}. By polynomializing nonlinearities into quadratic observables, the method preserves expressive power while keeping computations feasible. This balance enables accurate system identification and prediction within a finite yet rich dynamical framework. This approach outperforms SINDy under certain conditions. For example, as shown in \cite{klus2025quadratic}, QENDy outperforms SINDy because it exactly captures the rational dynamics by lifting the system into a higher-dimensional space, achieving machine-precision recovery. SINDy, on the other hand, fails unless the correct rational terms are manually added to the dictionary, leading only to approximate dynamics with higher errors. Nevertheless, QENDy also requires trajectory data along with corresponding time derivatives $\left(x_{i}, \dot{x}_{i}\right)_{i=1}^{m}$, where $m$ denotes the number of data points. Since numerical differentiation at the given points is highly sensitive to noise, the method may suffer from instability and reduced accuracy.
\par
In this manuscript, an integral form of QENDy, which we call iQENDy, is developed that eliminates the need for numerical differentiation of the data points. This modification is inspired by the work \cite{messenger2021weaksindy}, where the authors introduced a weak formulation to system identification. This work is shown to be more robust against noisy data than the original algorithm. See also \cite{rosenfeld2024occupation}.
\par
The remainder of this article is structured as follows. First, the original QENDy algorithm is reviewed, and the necessary tools required for its implementation are described. Next, the integral formulation of QENDy is introduced, highlighting the modifications made to enhance its performance. Finally, numerical experiments are presented to compare the effectiveness of both approaches, with evaluations carried out under noise-free and noisy conditions.

\section{QENDy Background}

Many data-driven techniques use the Koopman operator or its generator to analyze and predict the behavior of the nonlinear dynamical system \cite{gonzalez2021kernel}, \cite{brunton2021koopman}, and \cite{klus2024networks}. The Koopman operator lifts the underlying dynamics to an infinite-dimensional space of basis functions, also called the space of observables. It is usually a Hilbert space. Eigenvalues and eigenfunctions of the Koopman operator are then approximated to learn the characteristics of the dynamics. While the Koopman operator offers a linear structure for analysis, it simultaneously involves an infinite-dimensional setting. Koopman operator theory and Dynamic Mode Decomposition (DMD) \cite{schmid2010dmd} provides powerful tools for analyzing nonlinear dynamical systems through a linear framework. They have been successfully applied in fluid dynamics to identify coherent flow structures and predict complex behaviors, as well as in power grids and mechanical systems for stability assessment and fault detection. Beyond engineering, these methods have also been utilized in biology, neuroscience, and finance, where they enable data-driven modeling, forecasting, and control of high-dimensional processes.

QENDy, on the other hand, also allows the quadratic representation of the dynamics to give a finite-dimensional representation. The QENDy algorithm aims to learn the governing equations of the quadratic embedding of a nonlinear, complex system in a higher-dimensional space of observables. Suppose there are $m$ data points ${x_{i}}$ obtained from a single trajectory. For each data point, the corresponding time derivatives ${\dot{x}_{i}}$ are required, which are obtained through numerical approximation. In addition, a set of observables, or a dictionary, is defined, for example, ${ f_1, f_2, \ldots, f_K }$. Let $f$= $(f_1,f_2,...,f_K)$. Setting $z^{(i)} = f(x^{(i)})$ and $\dot{z}^{(i)} = J(x^{(i)})\dot{x}^{(i)}$, where $J$ denotes the Jacobian matrix of the observables, the objective is to determine the matrices $A$, $B$, and $C$ that minimize the prescribed loss function
\begin{equation}
L(A,B,C) = \sum_{i =1}^{m}\left\| \dot{z}^{\left ( i \right )}- A\left ( z^{\left ( i \right )}\otimes z^{\left ( i \right )}  \right ) - Bz^{\left ( i \right )}- C\right\|_{2}^{2}.
\end{equation}

Defining the data matrices
\begin{equation*}
\begin{split}
 Z_1 &= \begin{bmatrix}
z^{\left ( 1 \right )} & z^{\left ( 2 \right )} & ... & z^{\left ( m \right )} 
\end{bmatrix} \in \mathbb{R}^{K \times m}, \\
 Z_2 &= \begin{bmatrix}
z^{\left ( 1 \right )}\otimes z^{\left ( 1 \right )} & z^{\left ( 2 \right )}\otimes z^{\left ( 2 \right )} & ... & z^{\left ( m \right )}\otimes z^{\left ( m \right )} \\
\end{bmatrix}\in \mathbb{R}^{K^{2} \times m}, \\
\dot{Z}_1 &= \begin{bmatrix}
\dot{z}^{\left ( 1 \right )} & \dot{z}^{\left ( 2 \right )} & ... & \dot{z}^{\left ( m \right )} 
\end{bmatrix} \in \mathbb{R}^{K \times m},
\end{split}
\end{equation*}
The loss function can be written as 
\begin{equation}
\label{Loss Function}
L\left ( A, B, C \right )=\left\| \dot{Z} - AZ_2 - BZ_1 -C  \mathrm{1}_{m}^{\intercal}\right\|_{F}^{2},
\end{equation}
where $\left\| .\right\|_{F}$ is the Frobenius norm and $\mathrm{1}_{m}$ $\in$ $\mathbb{R}^{m}$ is the vector of ones. The optimal matrices A, B, and C are obtained through regression.
\par 
The QENDy algorithm determines the quadratic embedding from which the governing equations for the nonlinear dynamical system can be extracted, assuming there is a matrix $ G $ such that $ x = G z $. QENDy requires the time data matrix $x$, and the corresponding time derivatives $dX$
\[
X = [x^{(1)} \ \cdots \ x^{(m)}], \quad 
\dot{X} = [\dot{x}^{(1)} \ \cdots \ \dot{x}^{(m)}],
\]
, a set of basis functions, and the gradient of the function matrix. The performance of the QENDy algorithm also depends on the choice of basis functions. For example, QENDy gives a strong performance for the Thomas attractor when a basis of trigonometric functions, $\{x, \sin{x},\cos{x}\}$, is employed, but using Hermite polynomials as the basis did not yield satisfactory results. Choosing the right set of functions for the given data is still an open question for algorithms like SINDy, QENDy, and EDMD. However, using numerical differentiation to obtain the matrix $dX$ is highly susceptible to noise, which can significantly amplify errors in the regression process.

\section{INTEGRAL QENDy}

In what follows, we introduce an integrated version of QENDy that does not require the derivative matrix $dX$ or the basis-function Jacobian. In this algorithm, the loss function \eqref{Loss Function} is integrated using a windowed approach across consecutive time steps, enabling more accurate, localized estimation of the integrated quantities. Note that choosing the size of the integral windows $h$ is crucial. Smaller integral windows help the model learn more accurately as it trains on multiple trajectories; on the other hand, larger windows lead to more smoothing.
\par
{\footnotesize
\begin{equation}
\begin{split}
    L_{t_i}\left ( A, B, C \right ) =  \left\| \int_{t_{i-h}}^{t_i}\dot{Z}\,dt -    A\int_{t_{i-h}}^{t_i}Z_2\,dt-   B\int_{t_{i-h}}^{t_i}Z_1\,dt -C\int_{t_{i-h}}^{t_i} dt\right\|_{F}^{2}, 
\end{split}
\end{equation}}
 where $i \in \left\{ 1,2,...,m \right\} $, and thus 
 {\small
\begin{equation}
\begin{split}
L_{t_i}(A, B, C) = \Bigg\| & Z\big(X(t_i)\big) - Z\big(X(t_{i-h})\big) 
- A \int_{t_{i-h}}^{t_i} Z_2 \, dt \\
& - B \int_{t_{i-h}}^{t_i} Z_1 \, dt 
- C \, (t_i - t_{i-h}) \Bigg\|_F^2.
\end{split}
\end{equation}}

The iQENDy algorithm leverages an integrated regression setup by stacking the governing equation over all time indices $i$, resulting in a formulation of the same size as the original regression problem. Once the matrices $A$, $B$, and $C$ are calculated, the final model is given by
\begin{equation}
 Z\big(t_i\big) = A \int_{t_{i-h}}^{t_i} Z_2 \, dt 
 +  B \int_{t_{i-h}}^{t_i} Z_1 \, dt 
+ C \, (t_i - t_{i-h}) + Z\big(t_{i-h}\big) 
\end{equation}

This integrated version follows the same solution procedure but benefits from improved numerical stability, especially in the presence of noise. By integrating the system dynamics rather than relying on numerical derivatives, this approach reduces error amplification and yields more reliable parameter estimates. As a result, the integrated QENDy formulation performs better on noisy datasets.

\section{Interpretation of the Robustness in the Integral formulation}
In this section, we will understand the mathematical concepts behind the improved robustness in the integral formulation iQENDy. As mentioned in \cite{doi:10.1137/1021044}, differentiation approximation of noisy observations is an ill-posed problem, also see \cite{article}.For example consider a continuously differentiable function $x_1(t)$ with derivative
\begin{equation*}
y_1(t) = \frac{d}{dt}x_1(t).
\end{equation*}
Define a perturbation function
\begin{equation*}
x_2(t) = x_1(t) + N\sin{(\omega t)},
\end{equation*}
where $N$ is arbitrarily small and  $\omega$ is arbitrarily large. Then the difference between $x_1$ and $x_2$ in sup norm satisfies 
\begin{equation*}
\left\| x_2(t) - x_1(t) \right\|_{2} = \left|{N}\right|,
\end{equation*}
However, their derivatives satisfy 
\begin{equation*}
y_2(t) = x_1'(t) + N\omega\cos{(\omega t)},
\end{equation*}
and therefore
\begin{equation*}
\left\| y_2(t) - y_1(t) \right\|_{2} = \left|{N}\omega\right|.
\end{equation*}
This demonstrates that differentiation is not a continuous operator on the space of continuous functions, and hence numerical differentiation from noisy measurements is an ill-posed problem. In particular, small high-frequency perturbations in the measured trajectory may produce arbitrarily large perturbations in its derivative, resulting in poor conditioning of the regression problem and unstable or divergent reconstructed trajectories.
\par
On the other hand, integration averages the  high-frequency noise. It acts as a smoothing functional on data matrices. Say we have the same problem,

\begin{equation*}
x_2(t) = x_1(t) + N \sin(\omega t).
\end{equation*}

Integrating over a finite window $[t-\tau, t]$, we obtain
\begin{equation*}
\int_{t-\tau}^{t} x_2(s)\, ds= \int_{t-\tau}^{t} x_1(s)\, ds + N \int_{t-\tau}^{t} \sin(\omega s)\, ds.
\end{equation*}

The perturbation term evaluates to
\begin{equation*}
N \int_{t-\tau}^{t} \sin(\omega s)\, ds
=\frac{N}{\omega}
\left[ \cos\big(\omega (t-\tau)\big) - \cos(\omega t)\right].
\end{equation*}

Therefore,
\begin{equation*}
\left\|
\int_{t-\tau}^{t} x_2(s)\, ds -\int_{t-\tau}^{t} x_1(s)\, ds
\right\|_{2} \le \frac{2|N|}{\omega}.
\end{equation*}

Although numerical integration is marginally more computationally intensive than numerical differentiation, the associated gain in robustness and stability makes it preferable in noisy system identification settings.
\par

\section{Numerical Results}
In this section, a comparative analysis between the original QENDy algorithm, the SINDy algorithm, and our integrated formulation is presented. To assess their performance, we consider the Thomas attractor and the Duffing oscillator as representative nonlinear dynamical systems. First, QENDy and iQENDy are compared under noiseless and noisy conditions. Both methods are applied to the same dataset and basis, and the resulting reconstructions are evaluated side by side. Then, to make a fair comparison among all three algorithms, a larger dictionary is used, and sparsity is introduced in QENDy and iQENDy. Experiments are conducted primarily by reconstructing the first three variables using all algorithms, with and without noise, allowing a clear qualitative assessment of accuracy, stability, and long-term predictive behavior.  This setup allows us to highlight the improvements introduced by the integrated approach while maintaining a consistent evaluation framework with the original method. \newline
Remark: Note that total variation filtering is applied on noisy data before computing the numerical derivatives for SINDy and QENDy. In contrast, no filtering is used in the implementation of iQENDy.

\subsection{The modified Thomas attractor}
The modified Thomas attractor is given by
\begin{equation}
\begin{split}
\dot{x}_1 &= \sin(x_2) - \alpha x_1 - \beta x_2 \cos(x_1), \\
\dot{x}_2 &= \sin(x_3) - \alpha x_2 - \beta x_3 \cos(x_2), \\
\dot{x}_3 &= \sin(x_1) - \alpha x_3 - \beta x_1 \cos(x_3).
\end{split}
\end{equation}

\subsubsection{QENDy vs. iQENDy}
In this study, we employ the same basis functions, $\left \{ x,\sin(x),\cos(x)\right\}$, and utilize an identical number of data points derived from a single trajectory over a fixed time interval for both the original QENDy algorithm and iQENDy. The key distinction between the two methods is that the original QENDy algorithm requires time derivatives, which are approximated using a central difference scheme to compute the derivative matrix. It should be noted that the computation of numerical derivatives forms part of our implementation and was not employed in the original work. In their approach, the authors used prior knowledge of the governing dynamics (the Thomas attractor) to directly obtain the derivative data matrix, thereby avoiding the need for numerical differentiation routines, whereas the iQENDy approach employs Simpson's rule, a higher-order numerical integration method.
For all simulations, 8,000 data points spanning the time interval from 0 to 60 units are generated. In the figures below, we compare reconstructions of the first three state variables, $\left\{x_1, x_2, x_3 \right\}$, obtained using the original QENDy algorithm and the updated iQENDy formulation, both with and without added noise. The value of $ \alpha $ and $ \beta $ is taken to be 0.2 and 0 (Thomas attractor), respectively.

 \begin{figure}[H]
     \centering
     \includegraphics[width=0.6\linewidth]{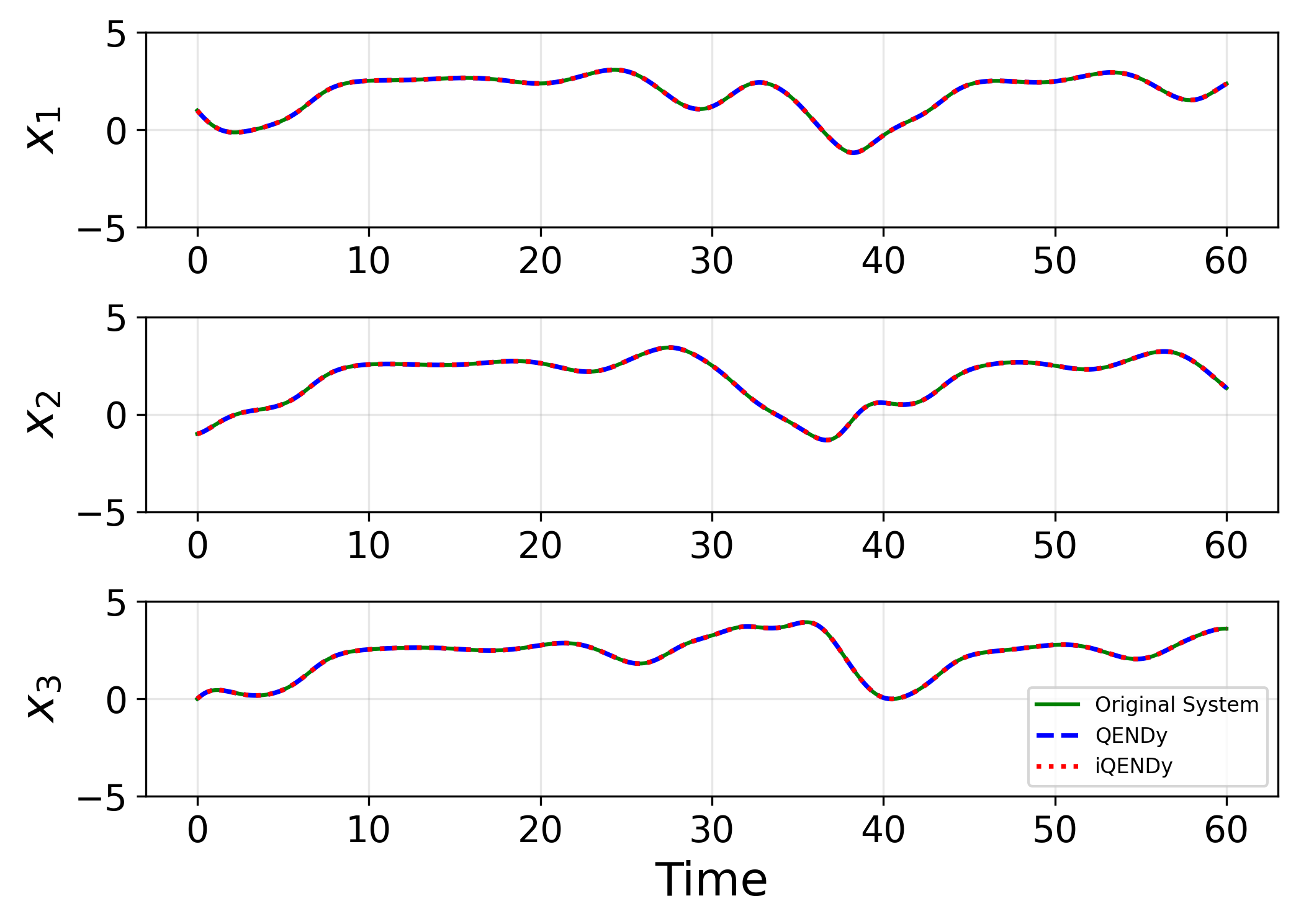}
     \caption{The above figure shows QENDy vs. iQENDy reconstruction of the first three state variables $(x_1,x_2,x_3)$ using noise-free data. Both models show good performance and perfectly fit the original trajectories. }
     \label{fig:placeholder}
 \end{figure}

\begin{figure}[H]
    \centering
    \includegraphics[width=0.6\linewidth]{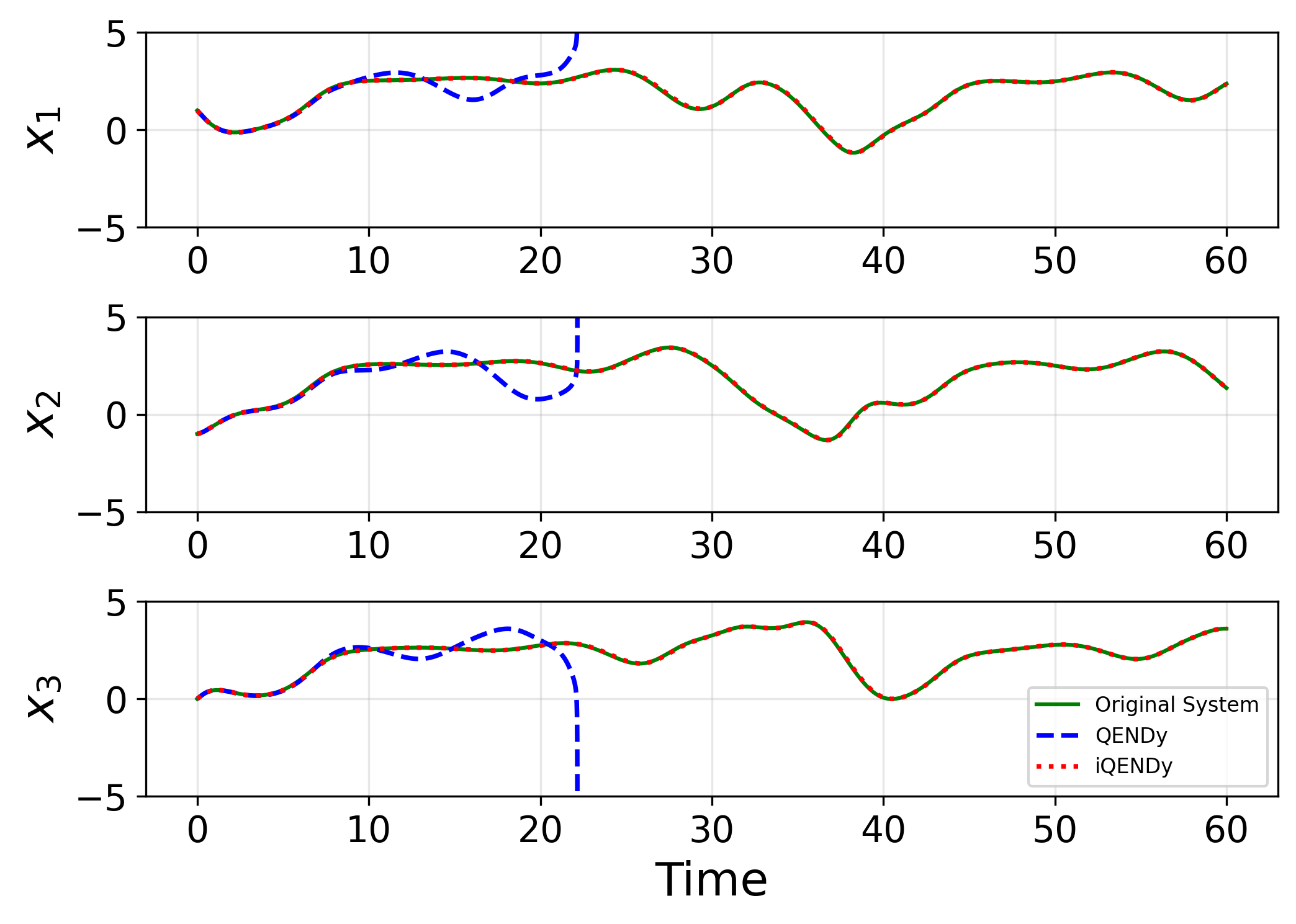}
    \caption{This figure shows reconstructions of the first three state variables using QENDy and iQENDy with small Gaussian noise added to the data matrix. It is clear that QENDy fails after a short period of time, but iQENDy still maintains accurate reconstruction.}
    \label{fig:placeholder}
\end{figure}

\subsubsection{SINDy vs. QENDy vs. iQENDy}Now, We employ a larger dictionary $\left \{ x,x^2,x^3\sin(x),\cos(x),x\sin(x),x\cos(x)\right\}$, and utilize the same number of data points derived from a single trajectory over
a fixed time interval for all three algorithms. Sparsity is introduced into the QENDy and iQENDy algorithms to facilitate better comparison with SINDy. Here, 10,000 data points
spanning the time interval from 0 to 60 units are generated. In
the figures presented below, we compare the reconstructions
of the first three state variables, $\left\{x_1, x_2, x_3 \right\}$, obtained using the original QENDy, SINDy, and iQENDy
formulation without noise and with added noise. The value
of $\alpha$ and $\beta$ is taken to be 0.2 and 0.004, respectively.

\begin{figure}[H]
    \centering
    \includegraphics[width=0.6\linewidth]{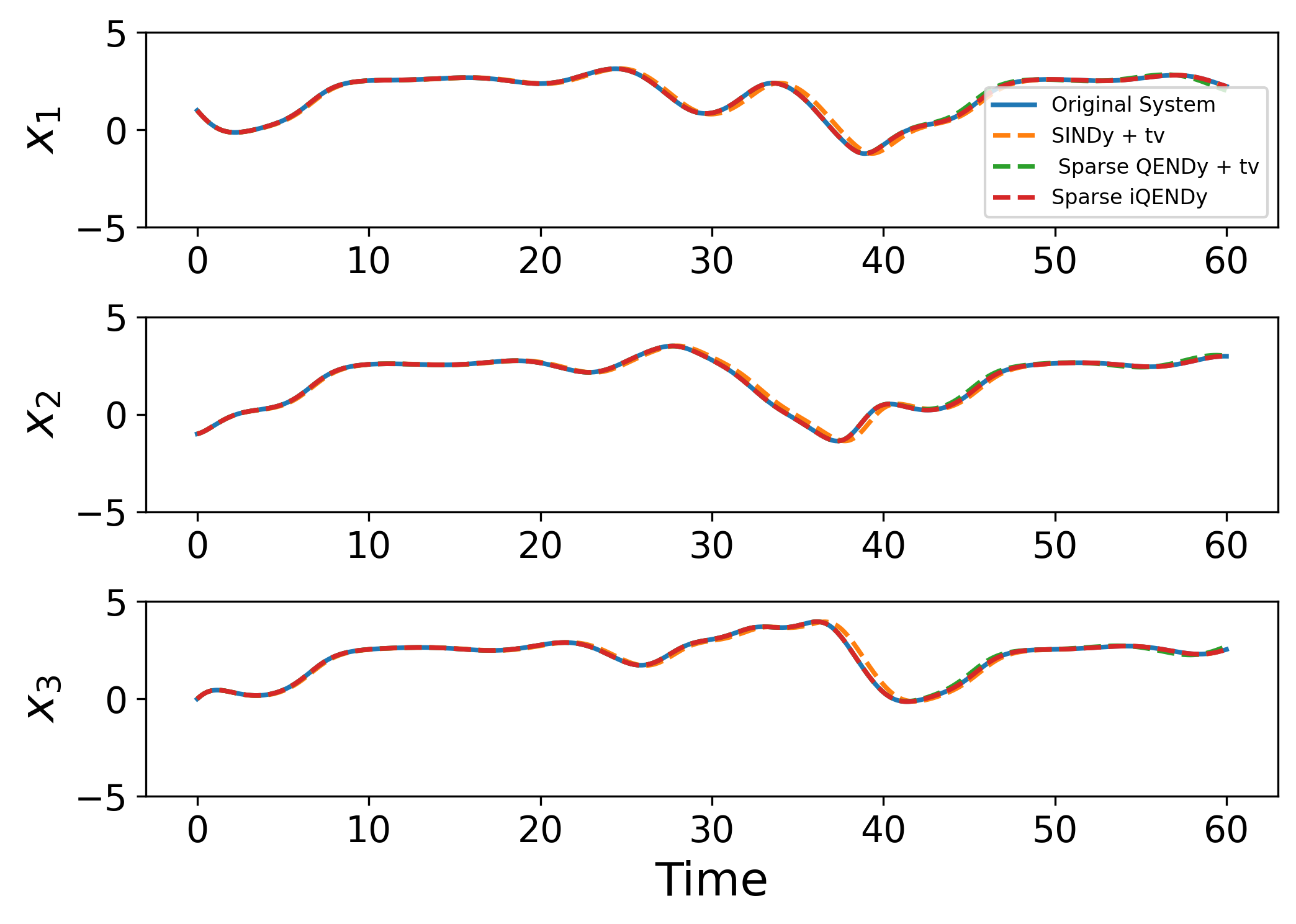}
    \caption{The reconstruction of the first three variables using sparse QENDy + tv (total variation filtering), sparse iQENDy +tv, and SINDy under clean conditions is depicted in this image. All three algorithms display good performance. }
    \label{fig:placeholder}
\end{figure}

\begin{figure}[H]
    \centering
    \includegraphics[width=0.6\linewidth]{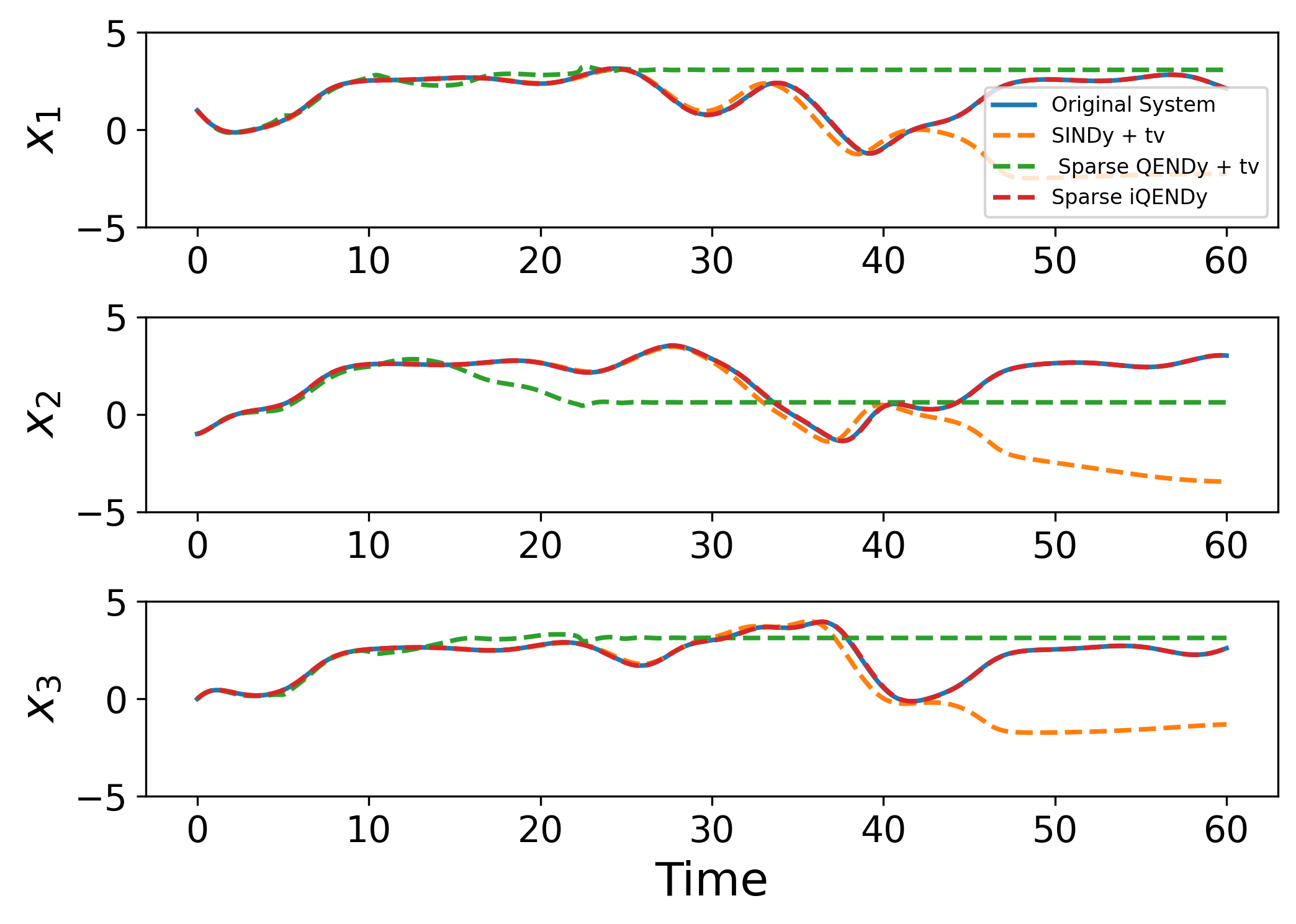}
    \caption{This figure shows a comparison of three algorithms under noisy conditions. It is evident that sparse iQENDy shows accurate performance over the longest interval.  }
    \label{fig:placeholder}
\end{figure}

\subsection{Duffing Oscillator}
The Duffing Oscillator system is given by 
\begin{equation}
\begin{split}
\dot{x}_1 &= x_2, \\
\dot{x}_2 &= -\delta x_2 - \alpha x_1 - \beta x_1^3,
\end{split}
\end{equation}
where $\delta$, $\alpha$, and $\beta$ are system parameters. 

\subsubsection{QENDy vs. iQENDy} To observe the chaotic behavior of the Duffing oscillator, the parameters are chosen as  
$\alpha = -1, \beta  = 2, \delta = 0.02$. From a trajectory with a prescribed initial condition, $5,000$ data points are extracted over a time interval of 0 to 60.  
For the observable space, a polynomial basis of the form  
$
\{\, x,\; \; x^{2},\;\; x^3\}
$
is employed. Both the QENDy and iQENDy algorithms are applied to reconstruct the first three state variables of the oscillator, and their performance is compared in the noise and noise-free settings.

\begin{figure}[H]
    \centering
    \includegraphics[width=0.6\linewidth]{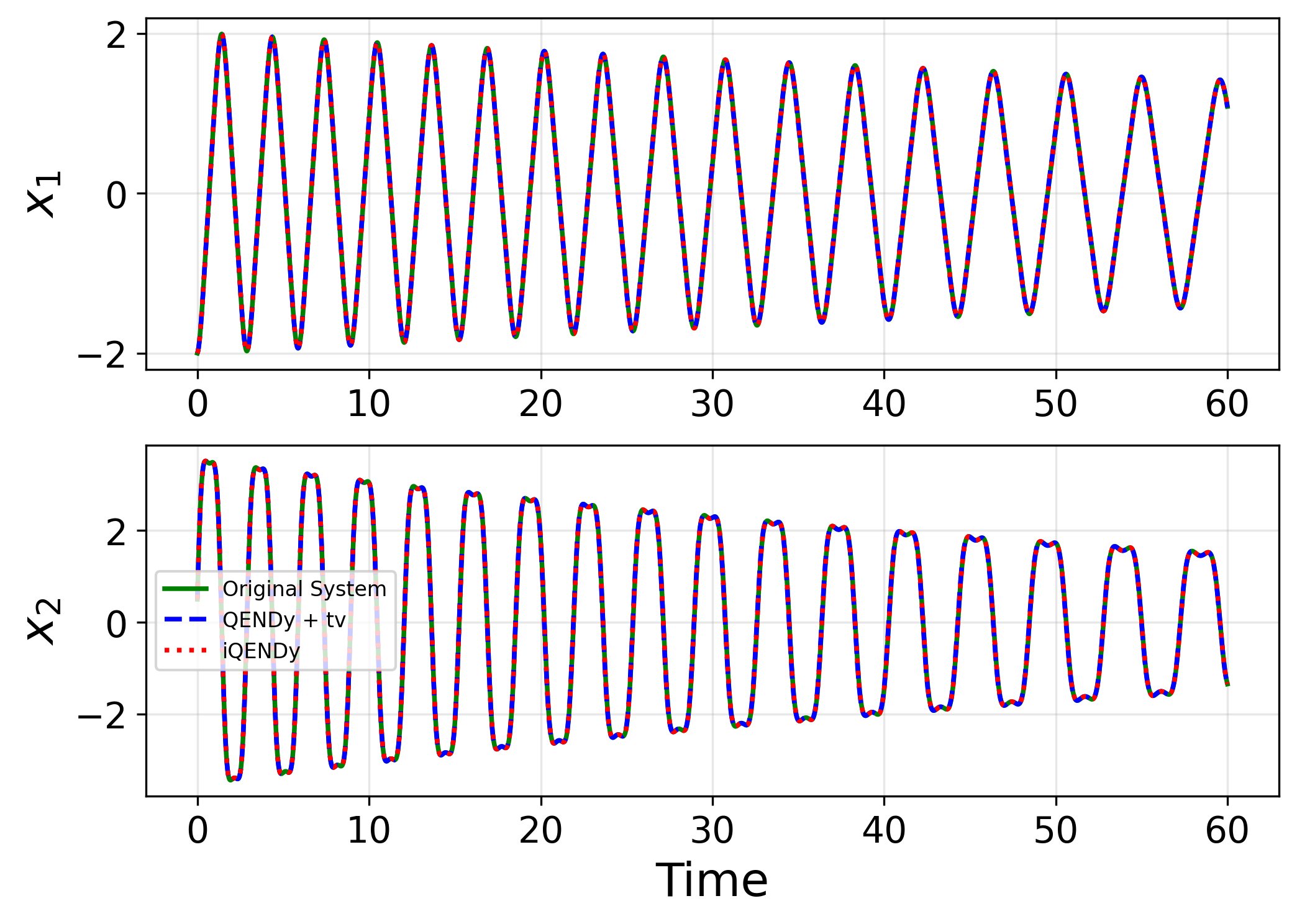}
    \caption{This figure shows the QENDy and iQENDy reconstructions of the first two state variables for the Duffing oscillator. When there is no noise, both models give accurate trajectories. }
    \label{fig:placeholder}
\end{figure}

\begin{figure}[H]
    \centering
    \includegraphics[width=0.6\linewidth]{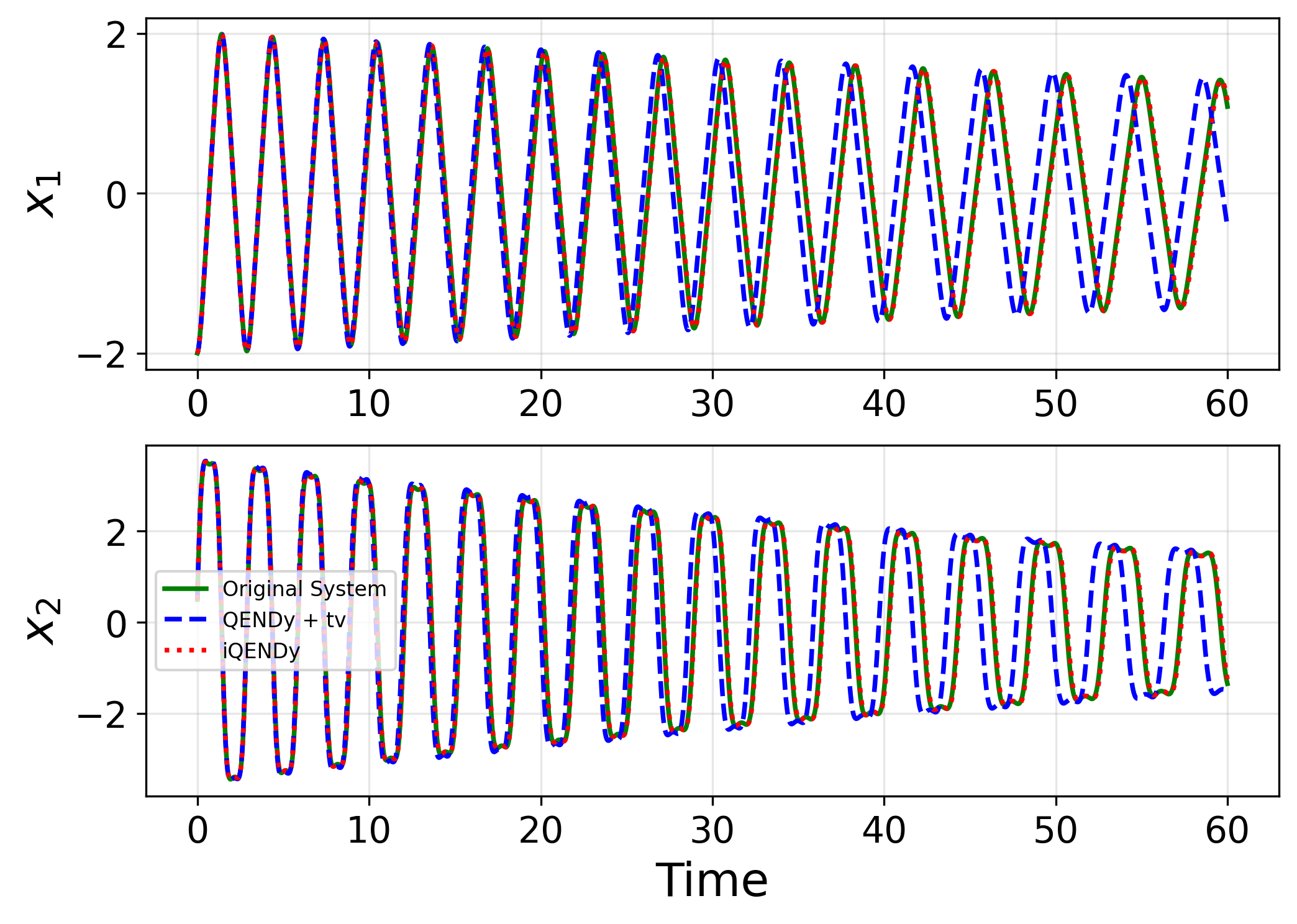}
    \caption{When a small amount of Gaussian noise is introduced, a phase shift in the QENDy reconstruction occurs after some time.}
    \label{fig:placeholder}
\end{figure}

\subsubsection{QENDy vs. SINDy vs. iQENDy}
Again, just like the Thomas attractor example, we use a larger dictionary  $\left \{ x,x^2,x^3\sin(x),\cos(x),x\sin(x),x\cos(x)\right\}$. 5000 data points are taken from a single trajectory over a time interval of 0 to 60. Below is the comparison between three algorithms with and without noise.
\begin{figure}[H]
    \centering
    \includegraphics[width=0.6\linewidth]{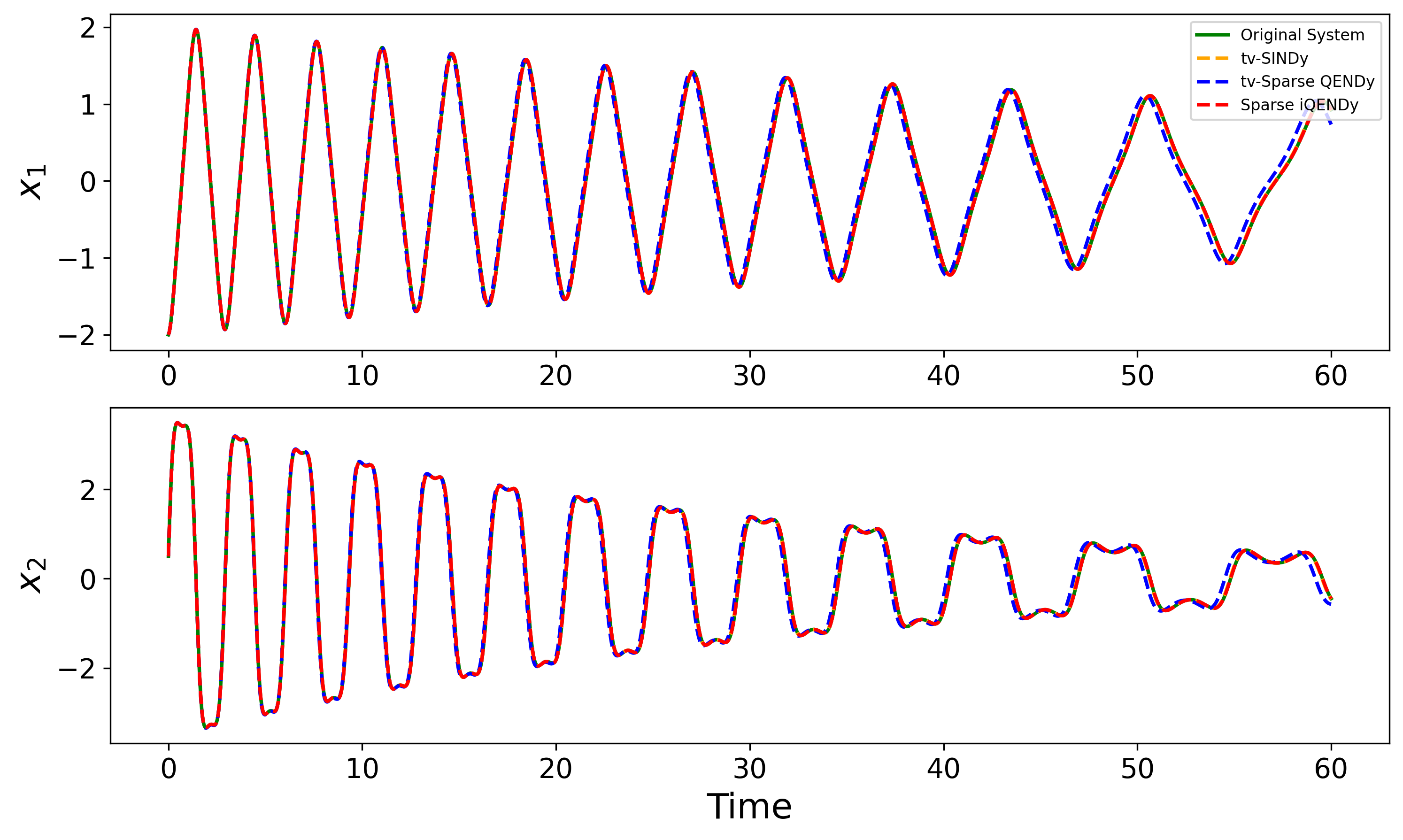}
    \caption{Reconstruction of the first two variables using sparse QENDy, sparse iQENDy, and SINDy is shown in the figure above. Even in noise-free conditions, a small phase shift is observed in the sparse QENDy reconstruction.  }
    \label{fig:placeholder}
\end{figure}
\begin{figure}[H]
    \centering
    \includegraphics[width=0.6\linewidth]{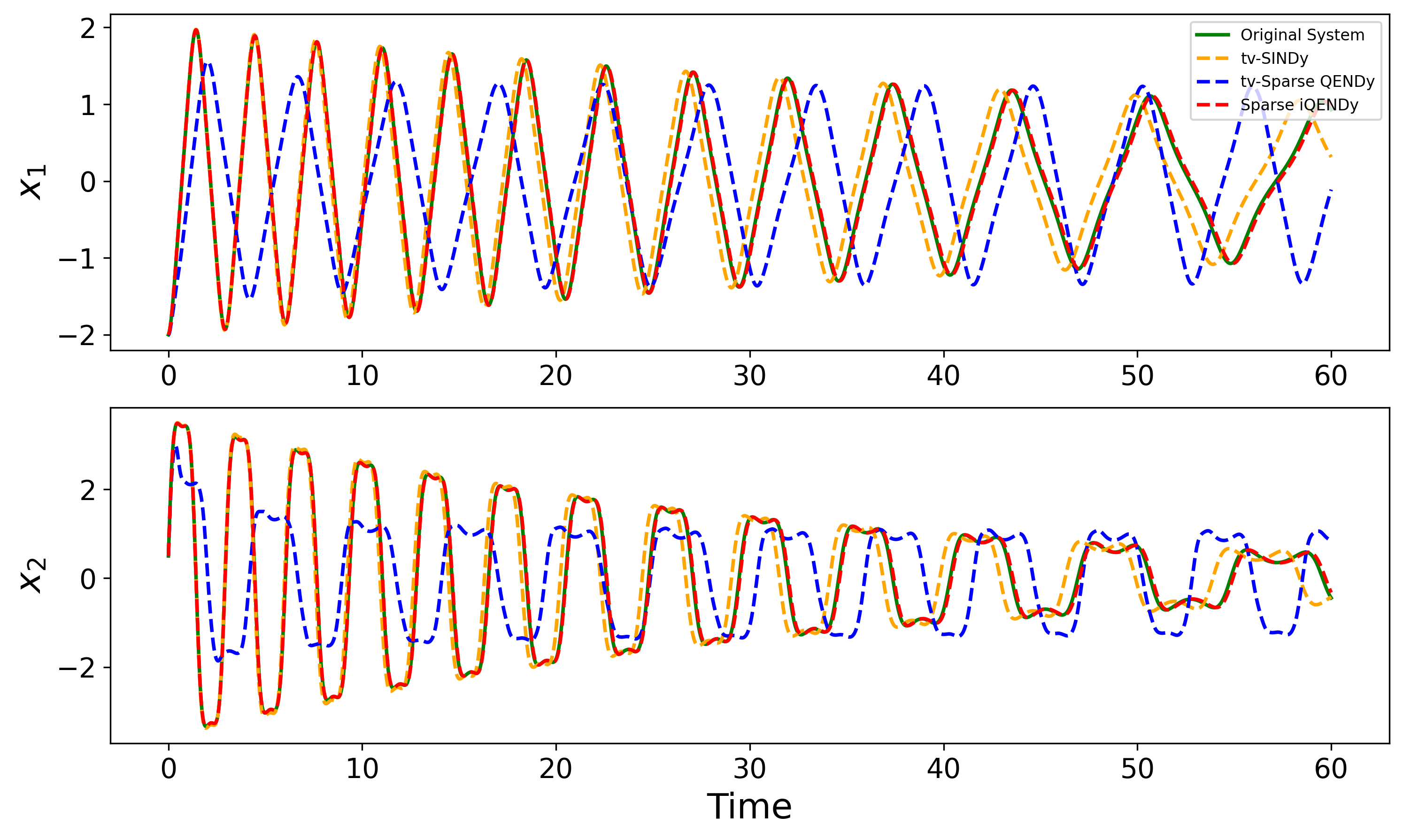}
    \caption{When noise is introduced into the data matrix, we observe phase shifts and poor performance from SINDy and sparse QENDy, whereas iQENDy remains stable.}
    \label{fig:placeholder}
\end{figure}
These results confirm that iQENDy outperforms QENDy and SINDy in the presence of noise. While QENDy and SINDy suffer from noise amplification during numerical differentiation, leading to phase drift and reduced accuracy, iQENDy remains stable. This highlights the effectiveness of the integral formulation in achieving robust system identification under noisy measurements.
\section{Generalization}
Learned models are better tested on a different initial condition than the one used for training to assess their capacity for generalization. Specifically, a new initial state is used to simulate the models generated by QENDy and iQENDy, and the resulting trajectories are compared with the system's actual dynamics. Instead of only fitting the training trajectory, this experiment evaluates whether the learned governing equations describe the underlying dynamics. The identified model provides an explicit quadratic polynomial description of the lifted dynamics for iQENDy. In the lifted coordinate space, this representation can be understood as a closed-form dynamical model that enables direct simulation for any initial conditions. The figure below shows the performance of QENDy and iQENDy on a test initial condition $x_1 = (1,2,1)$.
\begin{figure}[H]
    \centering
    \includegraphics[width=0.6\linewidth]{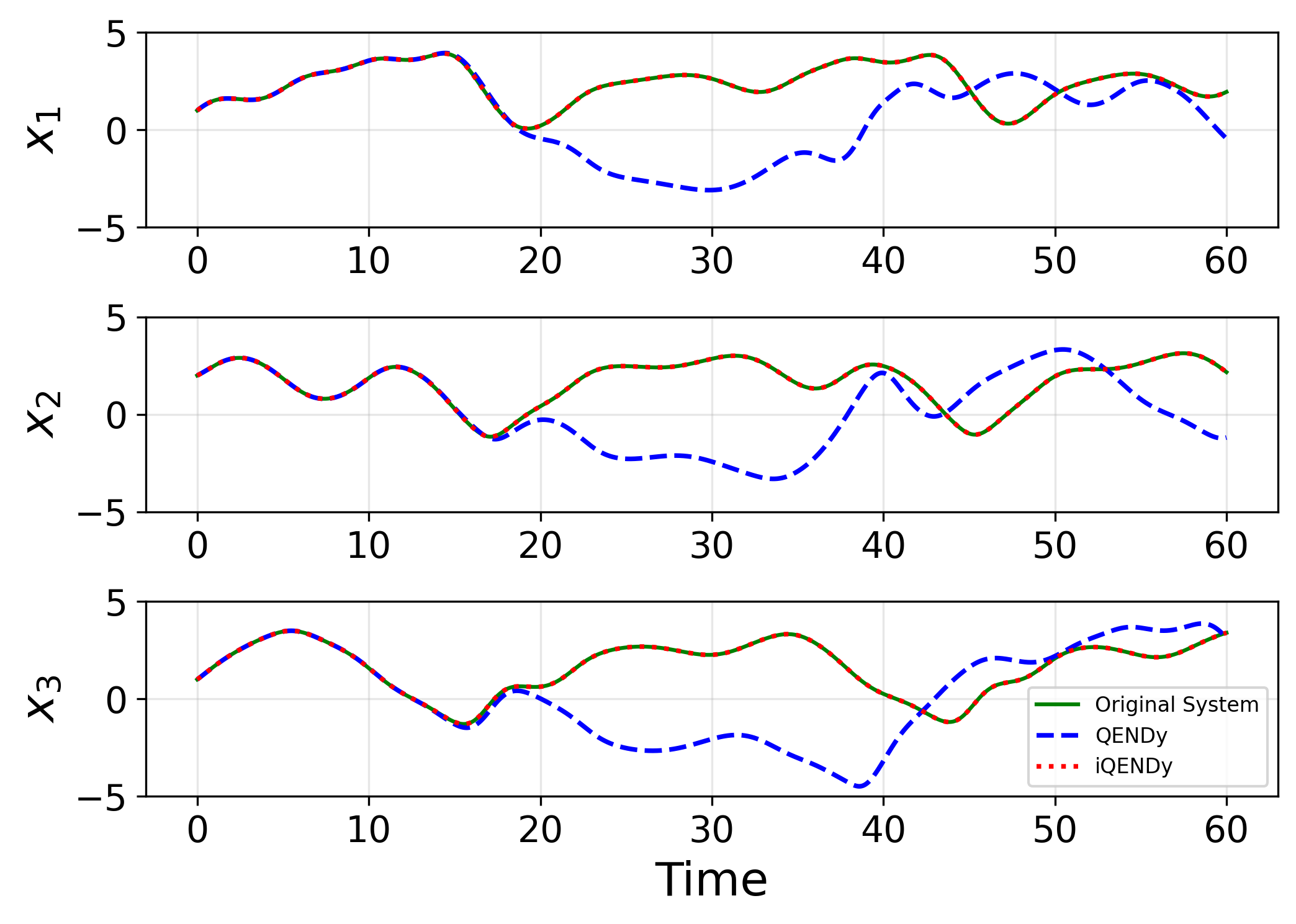}
    \caption{The above figure shows the QENDy and iQENDy reconstruction of the first three variables of the trajectory with the test initial condition $x_1 = (1,2,1)$. Both models were trained on initial condition $x_0 = (1,-1,0)$  }
    
\end{figure}

\section{Error Analysis}
To test the accuracy and robustness of the iQENDy algorithm, the reconstruction mean square error is compared for the Thomas attractor under increasing noise for each algorithm. These simulations have been run on the same time grid and with the same initial conditions to ensure a fair test. Below are two tables: one for QENDy vs. iQENDy and the other for sparse QENDy vs. iQENDy vs. SINDy.

\begin{table}[h]
\centering
\begin{tabular}{ccc}
\toprule
\textbf{Noise (\%)} & \textbf{QENDy-tv} & \textbf{iQENDy} \\
\midrule
0  & 2.49e-3 & 1.77e-10 \\
0.1  & 1.82e+1 & 2.99e-1 \\
1  & 8.07e+1 & 2.99e+0 \\
2 & 111.5 & 5.99 \\
5 & 178.5 & 14.9 \\
\bottomrule
\end{tabular}
\caption{\normalfont Mean squared error comparison between QENDy and iQENDy for the Thomas attractor under different noise levels.}
\label{tab:qendy_iqendy_error}
\end{table}

\begin{table}[H]
\centering
\begin{tabular}{cccc}
\toprule
\textbf{Noise (\%)} & \textbf{SINDy-tv} & \textbf{Sp. QENDy-tv } & \textbf{Sp. iQENDy } \\
\midrule
0  & 1.92e-1 & 2.15e-2 & 2.18e-4 \\
0.1  & 11.6 & 157.6 & 4.5 \\
1  & 80.3 & 1380.2 & 48.2 \\
2  & 140.9 & 2527.7 & 96.03 \\

\bottomrule
\end{tabular}
\caption{\normalfont Mean squared error comparison between SINDy, Sparse QENDy, and Sparse iQENDy for the Thomas attractor under different noise levels.}
\label{tab:sparse_methods_error}
\end{table}
A similar framework of iQENDy reconstruction error analysis can be applied to the Duffing oscillator system. The iQENDy method can be further explored by applying it to a variety of dynamical systems, with different values of parameters, including those with nonlinear, rational, or chaotic behavior, while also testing its performance under limited data availability.  At the same time, selecting an appropriate basis remains a crucial challenge when dealing with unseen dynamics, as an inadequate dictionary may limit the method’s ability to generalize. Adding dictionary learning tools appears to be a promising direction for future work, see \cite{GOYAL2024134158}.
\section{Conclusion}

In this work, we introduced an integrated formulation of the QENDy algorithm tailored for nonlinear dynamical systems and demonstrated its effectiveness on the Thomas system and chaotic Duffing oscillator examples. Through a detailed numerical comparison with the original QENDy approach, our method demonstrated improved accuracy and robustness in reconstructing system trajectories over a longer period. The integrated approach benefits from not using numerical derivatives, allowing it to better capture the nonlinear interactions inherent in complex dynamics and resulting in increased robustness. These results highlight the potential of the iQENDy formulation as a powerful tool for data-driven system identification.

\bibliographystyle{IEEEtran}
\bibliography{References}

\end{document}